\documentclass[a4paper,11pt,reqno]{amsart}

\usepackage{dsfont,amssymb}

\usepackage{hyperref}

\numberwithin{equation}{section}

\newtheorem{theorem}{Theorem}[section]
\newtheorem{corollary}[theorem]{Corollary}

\newtheorem{proposition}[theorem]{Proposition}
\theoremstyle{remark}

\newtheorem{example}[theorem]{Example}
\theoremstyle{definition}

\newcommand\bp{\begin{proof}}
\newcommand\ep{\end{proof}}
\newcommand\ee{\nopagebreak\mbox{\ }\hfill$\diamondsuit$}

\newcommand\Ad{\operatorname{Ad}}
\newcommand{\R}{{\mathbb R}}
\newcommand\T{{\mathbb T}}

\newcommand{\A}{{\mathcal A}}
\newcommand{\F}{{\mathcal F}}

\begin{document}

\title{Smooth crossed products of Rieffel's deformations}

\author[S. Neshveyev]{Sergey Neshveyev}

\email{sergeyn@math.uio.no}

\address{Department of Mathematics, University of Oslo,
P.O. Box 1053 Blindern, NO-0316 Oslo, Norway}

\thanks{The research leading to these results has received funding from the European Research Council
under the European Union's Seventh Framework Programme (FP/2007-2013) / ERC Grant Agreement no. 307663
}

\date{July 8, 2013; minor changes November 14, 2013}

\begin{abstract}
Assume $\A$ is a Fr\'echet algebra equipped with a smooth isometric action of a vector group~$V$, and consider Rieffel's deformation $\A_J$ of $\A$. We construct an explicit isomorphism between the smooth crossed products $V\ltimes\A_J$ and $V\ltimes\A$. When combined with the Elliott-Natsume-Nest isomorphism, this immediately implies that the periodic cyclic cohomology is invariant under deformation. Specializing to the case of smooth subalgebras of C$^*$-algebras, we also get a simple proof of equivalence of Rieffel's and Kasprzak's approaches to deformation.
\end{abstract}

\maketitle

\bigskip

\section*{Introduction}

The main goal of this note is to give a short proof of invariance of periodic cyclic cohomology under Rieffel's deformations. Particular cases of this result are, of course, well-known. For the noncommutative $2$-tori this was already shown by Connes in the foundational paper~\cite{C}. The result was extended to the higher dimensional noncommutative tori by Nest~\cite{Nest}. More recently, similar results have been obtained for $\theta$-deformations by Yamashita~\cite{Yam1} and Sangha~\cite{San}. A possible systematic way of approaching the question of invariance of periodic cyclic (co)homology is by using the Gauss-Manin connection, see e.g.~\cite{Yam2,Yas}, but in the analytic setting this usually involves significant technical difficulties. It is often more efficient to use crossed product decompositions.

Given a Fr\'echet algebra $\A$ with a smooth isometric action of a vector group~$V$, for Rieffel's deformation $\A_J$ of $\A$ we construct an isomorphism between the smooth crossed products $V\ltimes\A_J$ and $V\ltimes\A$. The existence of such an isomorphism on the C$^*$-algebra level is known~\cite{Kas,BNS}, but the proof of this existence has been rather indirect and relied heavily on the C$^*$-algebra technique. As it turns out, the origin of this isomorphism could not be easier: both smooth crossed products are naturally represented on the space $S(V;\A)$ of $\A$-valued Schwartz functions on~$V$, and their images under these representations coincide.

The isomorphism $V\ltimes\A_J\cong V\ltimes\A$ gives an embedding of~$\A_J$ into the multiplier algebra of~$V\ltimes\A$. For $\theta$-deformations a formula for such an embedding in terms of the decomposition of $\A$ into spectral subspaces is easy to guess, which was already used in the work of Connes and Landi~\cite{CL}. For general Rieffel's deformations, when there are no nonzero spectral subspaces, it is impossible to write down such a formula, yet the isomorphism $V\ltimes\A_J\cong V\ltimes\A$ has an explicit and relatively simple form.

In the second part of this note we consider smooth subalgebras $\A\subset A$ of C$^*$-algebras. For C$^*$-algebras, a different approach to deformation has been proposed by Kasprzak~\cite{Kas}. In his approach the existence of an isomorphism $V\ltimes A_J\cong V\ltimes A$ is taken as part of the definition of $A_J$, so that~$A_J$ is from the beginning defined as a subalgebra of  $M(V\ltimes A)$. Concretely, elements of $A_J$ can be obtained using either Landstad's theory~\cite{Kas} or certain quantization maps $A\to M(V\ltimes A)$~\cite{BNS,NTdef2}. Equivalence of two approaches has been proved in~\cite{BNS}, but the proof was far from straightforward. Using our explicit isomorphism $V\ltimes\A_J\cong V\ltimes\A$ we can now give a very simple proof. To complete the picture, we also describe the quantization maps $A\to M(V\ltimes A)$ in Rieffel's setting.

\bigskip

\section{Fr\'echet algebras}

Assume $\A$ is a Fr\'echet algebra with a smooth isometric action $\alpha$ of a vector group $V\cong\R^d$. Following Rieffel~\cite{Ri1}, by saying that the action is isometric we mean that the topology on $\A$ is defined by a sequence of $\alpha$-invariant seminorms. The assumption of smoothness means that for every $a\in\A$ the function $V\to\A$, $x\mapsto\alpha_x(a)$, is differentiable and the partial differentiation in the direction $X\in V$ at $x=0$ defines a bounded operator on $\A$ for every $X$.

Fix a scalar product $\langle\cdot,\cdot\rangle$ on $V$. Consider the space $S(V;\A)$ of $\A$-valued Schwartz functions on~$V$. It can be made into a Fr\'echet algebra by defining the convolution product by
$$
(f*g)(x)=\int_Vf(y)\alpha_y(f(x-y))dy.
$$
Following \cite{ENN} we denote the space  $S(V;\A)$ with this product by $V\ltimes_\alpha\A$ and call it the smooth crossed product of $\A$ by the action of $V$.

Let $J$ be a skew-symmetric operator on $V$. Then Rieffel's deformation $\A_J$ of $\A$ is the Fr\'echet space $\A$ equipped with the new product
$$
a\times_Jb=\int_{V\times V}\alpha_{Jx}(a)\alpha_y(b)e(x\cdot y)dx\,dy,
$$
where $e(x\cdot y)$ stands for $e^{2\pi i\langle x,y\rangle}$ and the integral is understood in the oscillatory sense~\cite{Ri1}. The automorphisms $\alpha_x$ of $\A$ remain automorphisms of $\A_J$ and define an action of $V$ on $\A_J$, which we denote by $\alpha^J$.

Define the Fourier transform $S(V;\A)\to S(V;\A)$ by
$$
\hat\xi(x)=\int_V \xi(y)e(-x\cdot y)dy.
$$
We have a representation $\pi$ of $\A$ on $S(V;\A)$ defined by
$$
(\pi(a)\xi)(x)=\alpha_{-x}(a)\xi(x).
$$
We also have a representation $\pi_J$ of $\A_J$ on $S(V;\A)$ defined by $\pi_J(a)\xi=\alpha(a)\times_J\xi$, where $\alpha(a)$ is the $\A$-valued function on $V$ given by $\alpha(a)(x)=\alpha_{-x}(a)$, and the deformed product $\times_J$ for $\A$-valued functions is defined using the action of $V$ on itself by left translations. In other words,
$$
(\pi_J(a)\xi)(x)=\int_{V\times V}\alpha_{Jy-x}(a)\xi(x-z)e(y\cdot z)dy\,dz
=\int_V\alpha_{Jy-x}(a)\hat\xi(y)e(x\cdot y)dy,
$$
where the second equality is justified by the computation before \cite[Proposition~3.1]{Ri1}. Note that the last integral, as well as all other integrals we will encounter from now on, is understood in the usual sense, we no longer need oscillation.

We also have a representation of $V$ on $S(V;\A)$ by the operators $\lambda_x\otimes1$ of left translation, so $$
((\lambda_x\otimes1)\xi)(y)=\xi(y-x).
$$
The representations $\pi$ of $\A$ and $\pi_J$ of $\A_J$ are covariant with respect to this representation, meaning that
$$
(\lambda_x\otimes1)\pi(a)(\lambda_{-x}\otimes1)=\pi(\alpha_x(a)),
$$
and similarly for $\pi_J$ and $\alpha^J$. These covariant representations define representations of the smooth crossed products $V\ltimes_\alpha \A$ and $V\ltimes_{\alpha^J} \A_J$ on $S(V;\A)$, which we continue to denote by $\pi$ and $\pi_J$, respectively. Since we are not dealing with isometric representations on Banach spaces, we have to check that the representations $\pi$ and $\pi_J$ are indeed well-defined, but this is clear from the following identities:
$$
(\pi(f)\xi)(x)=\int_V\alpha_{-x}(f(y))\xi(x-y)dy
$$
and
\begin{align}
(\pi_J(f)\xi)(x)&=\int_V(\pi_J(f(z))(\lambda_z\otimes1)\xi)(x)dz\nonumber\\
&=\int_{V\times V}\alpha_{Jy-x}(f(z))[{(\lambda_z\otimes1)\xi}]^\wedge(y)e(x\cdot y)dy\,dz\nonumber\\
&=\int_{V\times V}\alpha_{Jy-x}(f(z))\hat\xi(y)e(-z\cdot y)e(x\cdot y)dy\,dz\nonumber\\
&=\int_{V}\alpha_{Jy-x}(\hat f(y))\hat\xi(y)e(x\cdot y)dy.\label{erep0}
\end{align}

Define an operator $\Theta_J$ on $S(V;\A)$ by
$$
\Theta_J(f)(x)=\int_V\alpha_{Jy}(\hat f(y))e(x\cdot y)dy.
$$
This operator is invertible, with inverse equal to $\Theta_{-J}$.

\begin{theorem}\label{tmain}
For every $f\in S(V;\A)$ we have $\pi_J(f)=\pi(\Theta_J(f))$. Hence $\Theta_J$ defines an isomorphism $V\ltimes_{\alpha^J} \A_J\cong V\ltimes_{\alpha} \A$.
\end{theorem}

\bp We compute:
\begin{align*}
(\pi(\Theta_J(f))\xi)(x)&=\int_V\alpha_{-x}(\Theta_J(f)(z))\xi(x-z)dz\\
&=\int_{V\times V}\alpha_{Jy-x}(\hat f(y))\xi(x-z)e(y\cdot z)dy\,dz\\
&=\int_{V\times V}\alpha_{Jy-x}(\hat f(y))\xi(z)e(-y\cdot z)e(x\cdot y)dy\,dz\\
&=\int_{V}\alpha_{Jy-x}(\hat f(y))\hat\xi(y)e(x\cdot y)dy.
\end{align*}
The last expression is exactly \eqref{erep0}, hence $\pi_J(f)=\pi(\Theta_J(f))$. If $\A$ is unital, then $\pi$ is injective, and we conclude that $\Theta_J$ is an isomorphism of the algebras $V\ltimes \A_J$ and $V\ltimes \A$. In general, we could extend, using the same formulas as before, the representations $\pi_J$ and $\pi$ to representations on the space $S(V;\A^\sim)\cong S(V;\A)\oplus S(V)$, where $\A^\sim$ is the unitization of $\A$. Then $\pi$ becomes injective and we still have $\pi_J(f)=\pi(\Theta_J(f))$, so we can  again conclude that $\Theta_J$ is an isomorphism of Fr\'echet algebras.
\ep

We remark that when $J$ is invertible, the operator $\Theta_J$ is the composition of the Fourier transform with an operator~$M$ introduced in \cite{BM}. But as we see, the origin of our operator is completely straightforward and its purpose is rather different from~\cite{BM}: the main result in \cite{BM} relates the stabilization of~$\A_J$ to a twisted crossed product of $\A$ by~$V$ in the C$^*$-algebraic setting.

\smallskip

Combining the above theorem with the Elliott-Natsume-Nest isomorphism~\cite{ENN} $$HP^*(V\ltimes\A)\cong HP^{*+d}(\A)$$ for the periodic cyclic cohomology with continuous cochains, we immediately get the following.

\begin{corollary}
We have $HP^*(\A_J)\cong HP^*(\A)$.
\end{corollary}

\bigskip

\section{\texorpdfstring{C$^*$}{C*}-algebras}

Assume now that $A$ is a C$^*$-algebra and $\alpha$ is a continuous action of $V\cong\R^d$ on $A$. Denote by $\A\subset A$ the algebra of smooth vectors for the action. Then the deformation $\A_J$ of $\A$ can be completed to a C$^*$-algebra $A_J$: the C$^*$-norm on $\A_J$ is defined by considering the representation $\pi_J$ as a representation on the right Hilbert $A$-module $L^2(V)\otimes A$, see~\cite{Ri1}.

\begin{theorem}
The isomorphism $\Theta_J$ of smooth crossed products extends to an isomorphism of the C$^*$-algebra crossed products $V\ltimes_{\alpha^J} A_J$ and $V\ltimes_\alpha A$.
\end{theorem}

\bp Since the representation $\pi_J$ of $A_J$ on the Hilbert $A$-module $L^2(V)\otimes A$ is covariant with respect to the left regular representation of $V$, it defines a representation of $V\ltimes A_J$. Therefore $\Theta_J$ extends to a homomorphism $V\ltimes A_J\to V\ltimes A$. Since $A=(A_J)_{-J}$, for the same reason the map~$\Theta_{-J}$ extends to a homomorphism $V\ltimes A\to V\ltimes A_J$. As the maps $\Theta_J$ and $\Theta_{-J}$ are inverse to each other on smooth crossed products, this gives the result.
\ep

A different approach to deformation of C$^*$-algebras has been developed in~\cite{Kas} and extended in~\cite{BNS,NTdef2}. It works for actions of arbitrary locally compact quantum groups and measurable cocycles on their duals. In order to describe it for $V$, it is convenient to start with a more general case of a locally compact abelian group $G$ and a continuous $2$-cocycle $\Omega\colon\hat G\times\hat G\to\T$. We will mainly follow the conventions in~\cite{NTdef2}, which are slightly different from those in~\cite{BNS}. It will be convenient though to write some of the formulas using $L^2(\hat G)$ instead of~$L^2(G)$.

Fix a Haar measure on $G$. Define the Fourier transform $\F\colon L^2(G)\to L^2(\hat G)$ by
$$
(\F f)(\chi)=\hat f(\chi)=\int_Gf(g)\overline{\chi(g)}dg.
$$
We normalize the Haar measure on $\hat G$ so that $\F$ becomes isometric. Consider the operators
$$
\lambda^\Omega_\chi=\lambda_\chi\overline{\Omega(\chi,\cdot)}\ \ \text{on}\ \ L^2(\hat G),
$$
where $\lambda_\chi$ are the operators of the left regular representation and $\overline{\Omega(\chi,\cdot)}$ is considered as the operator of multiplication by the function $\overline{\Omega(\chi,\cdot)}$. We have $\lambda^\Omega_{\chi\chi'}=\Omega(\chi,\chi')\lambda^\Omega_{\chi}\lambda^\Omega_{\chi'}$. Denote by $C^*_r(\hat G;\Omega)$ the C$^*$-algebra obtained as the norm closure of the space of operators of the form $\int_{\hat G}f(\chi)\lambda^\Omega_\chi d\chi$ with $f\in L^1(\hat G)$, and denote by $W^*(\hat G;\Omega)\subset B(L^2(\hat G))$ the von Neumann algebra it generates. For every normal linear functional $\nu\in W^*(\hat G;\bar\Omega)_*$ define a ``quantization map''
$$
T_\nu\colon C_0(G)\to C^*_r(\hat G;\Omega)\ \ \text{by}\ \
T_\nu(f)=(\iota\otimes\nu)(\hat W\Omega(\F f\F^*\otimes1)(\hat W\Omega)^*),
$$
where $\hat W$ is the multiplicative unitary of $\hat G$, so $(\hat W\xi)(\chi,\chi')=\xi(\chi,\chi^{-1}\chi')$. Explicitly, if $f$ lies in the Fourier algebra $A(G)\subset C_0(G)$ of $G$, so it is the inverse Fourier transform of a function $\hat f\in L^1(\hat G)$, then
$$
T_\nu(f)=\int_{\hat G}\hat f(\chi)\nu(\lambda^{\bar\Omega}_\chi)\lambda^{\Omega}_\chi d\chi.
$$
This follows from the identity $\hat W\Omega(\lambda_\chi\otimes1)(\hat W\Omega)^*=\lambda^\Omega_\chi\otimes\lambda^{\bar\Omega}_\chi$.
For every C$^*$-algebra $A$ the map $T_\nu\otimes\iota$ extends to a well-defined map $$M(C_0(G)\otimes A)\to M(C^*_r(\hat G;\Omega)\otimes A)$$
that is strictly continuous on the unit ball.

Assume now that we are given a continuous action $\alpha$ of $G$ on a C$^*$-algebra $A$. We view it as a homomorphism $\alpha\colon A\to M(C_0(G)\otimes A)$, so that $\alpha(a)(g)=\alpha_{-g}(a)$. We then define the deformation
$$
A_\Omega\subset M(C^*_r(\hat G;\Omega)\otimes A)
$$
as the C$^*$-algebra generated by the elements $(T_\nu\otimes\iota)\alpha(a)$ for all $a\in A$ and $\nu\in W^*(\hat G;\bar\Omega)_*$. This C$^*$-algebra carries a continuous action $\alpha^\Omega$ of $G$ defined by $\alpha^\Omega_g=\Ad(\F\lambda_g\F^*\otimes1)$, and the maps $(T_\nu\otimes\iota)\alpha \colon A\to A_\Omega$ are $G$-equivariant. Using the representation $g\mapsto \F\lambda_g\F^*\otimes1$ of $G$, we get a representation of the crossed product $G\ltimes_{\alpha^\Omega}A_\Omega$ on the right Hilbert $A$-module $L^2(\hat G)\otimes A$. Then \cite[Theorem~3.1]{BNS}, or \cite[Theorem~3.9]{NTdef2}, can be formulated by saying that this representation is faithful, and the map $\Ad(\F^*\otimes1)$ defines an isomorphism $G\ltimes_{\alpha^\Omega}A_\Omega\cong G\ltimes_\alpha A$,
if we identify $G\ltimes A$ with the norm closure of $\alpha(A)(C^*_r(G)\otimes1)$, that is, with the image of $G\ltimes A$ under its standard representation on the Hilbert $A$-module $L^2(G)\otimes A$.

It is easy to check that if we identify $G\ltimes_{\alpha^\Omega}A_\Omega$ with its image in the algebra of operators on the Hilbert $A$-module $L^2(\hat G)\otimes A$, then the dual action of $\hat G$ is defined by the automorphisms
$
\Ad(\Omega(\cdot,\chi)\lambda^*_\chi\otimes1).
$
Using the isomorphism $G\ltimes_{\alpha^\Omega}A_\Omega\cong G\ltimes_\alpha A$ we then get an action on $G\ltimes_\alpha A$, which we call the twisted dual action and denote by $\hat\alpha^\Omega$. Thus,
$$
\hat\alpha^\Omega_\chi=\Ad(\F^*\Omega(\cdot,\chi)\lambda^*_\chi\F\otimes1)
=\Ad(\F^*\Omega(\cdot,\chi)\F\otimes1)\circ\hat\alpha_\chi.
$$
Then an alternative description of $A_\Omega$ (or rather of $(\F^*\otimes1)A_\Omega(\F\otimes1)$), which is the original definition of Kasprzak~\cite{Kas} modulo replacing $\Omega$ by $\bar\Omega$, is that this is a unique $G$-invariant C$^*$-subalgebra of $M(G\ltimes A)^{\hat\alpha^\Omega}$ such that together with the embedding $C^*(G)\subset M(G\rtimes A)$ we get a decomposition $G\ltimes A=G\ltimes A_\Omega$ with respect to which the action $\hat\alpha_\Omega$ becomes the dual action on $G\ltimes A_\Omega$. The elements of this C$^*$-algebra can be abstractly characterized by Landstad's conditions, see~\cite{Kas} for details.

Finally, observe that the twisted dual action takes a simple form when $\Omega$ is a bi-character rather than just a $2$-cocycle. Indeed, in this case there exists a continuous homomorphism $r_\Omega\colon\hat G\to G$ such that
$$
\chi'(r_\Omega(\chi))=\Omega(\chi',\chi)\ \ \text{for all}\ \ \chi,\chi'\in\hat G.
$$
Then $\Omega(\cdot,\chi)$ is the operator of multiplication by the character $r_\Omega(\chi)$, so $\F^*\Omega(\cdot,\chi)\F=\lambda_{-r_\Omega(\chi)}$, and we get
\begin{equation}\label{edual}
\hat\alpha^\Omega_\chi=\Ad(\lambda_{-r_\Omega(\chi)}\otimes1)\circ\hat\alpha_\chi.
\end{equation}

\medskip

In order to illustrate the above definitions, consider a simple example.

\begin{example}\label{exap}
Assume that the action $\alpha$ is almost periodic, so that the spectral subspaces $$A_\chi=\{a\in A\mid \alpha_g(a)=\overline{\chi(g)}a\ \text{for all}\  g\in G\}$$ span a dense subalgebra of $A$. Then $A_\Omega\subset M(C^*_r(\hat G;\Omega)\otimes A)$ is the closed linear span of elements of the form $\lambda^\Omega_\chi\otimes a$ for $a\in A_\chi$ and $\chi\in\hat G$, since $(T_\nu\otimes\iota)\alpha(a)=\nu(\lambda^{\bar\Omega}_\chi)(\lambda^\Omega_\chi\otimes a)$ for $a\in A_\chi$, and we have
$$
(\lambda^\Omega_\chi\otimes a)(\lambda^\Omega_{\chi'}\otimes a')=\overline{\Omega(\chi,\chi')}
(\lambda^\Omega_{\chi\chi'}\otimes aa')\ \ \text{for}\ \ a\in A_\chi,\ a'\in A_{\chi'}.
$$
If in addition $\Omega$ is a skew-symmetric bi-character, so that $\overline{\Omega(\chi,\cdot)}=r_\Omega(\chi)$, then the embedding $\Ad(\F^*\otimes 1)$ of $A_\Omega$ into the multiplier algebra of $G\ltimes A=\overline{\alpha(A)(C^*_r(G)\otimes 1)}$ is given by
$$
\lambda^\Omega_\chi\otimes a\mapsto \F^*\lambda_\chi\overline{\Omega(\chi,\cdot)}\F\otimes a=\alpha(a)(\lambda_{-r_\Omega(\chi)}\otimes1)=(\lambda_{-r_\Omega(\chi)}\otimes1)\alpha(a)\ \ \text{for}\ \ a\in A_\chi.
$$
For $G=\T^n$ this is exactly the embedding used by Connes and Landi~\cite{CL} to construct a representation of the $\theta$-deformation of $A$ from a covariant representation of $A$.
\ee
\end{example}

Return to the case $G=V\cong\R^d$. We identify $\hat V$ with $V$ using the pairing $e(x\cdot y)$ and define a $2$-cocycle $\Omega_J$ on $\hat V=V$ by
$
\Omega_J(x,y)=e(x\cdot Jy).
$

\begin{theorem}\label{tdual}
The map $\Theta_J\colon M(V\ltimes_{\alpha^J} A_J)\to M(V\ltimes_\alpha A)$ defines an isomorphism of the C$^*$-algebras $A_J\subset M(V\ltimes_{\alpha^J} A_J)$ and $(\F^*\otimes 1)A_{\Omega_J}(\F\otimes1)\subset M(V\ltimes_\alpha A)$.
\end{theorem}

By the preceding discussion, here we identify $V\ltimes A$ with its image under the standard representation $\pi$ of $V\ltimes A$ on the Hilbert module $L^2(V)\otimes A$. Since $\pi\Theta_J=\pi_J$ by Theorem~\ref{tmain}, an equivalent way of formulating the above theorem is by saying that the map $a\mapsto \Ad(\F\otimes 1)\pi_J(a)$ defines an isomorphism $A_J\cong A_{\Omega_J}$.

\bp[Proof of Theorem~\ref{tdual}] It is enough to show that the isomorphism $\Theta_J\colon V\ltimes A_J\to V\ltimes A$ intertwines the dual action $\widehat{\alpha^J}$ on $V\ltimes A_J$ with the twisted dual action $\hat\alpha^{\Omega_J}$ on $V\ltimes A$. For $f\in S(V;\A)\subset V\ltimes A_J$ the dual action $\widehat{\alpha^J}$ is defined by
$$
\widehat{\alpha^J_y}(f)(x)=e(-x\cdot y)f(x).
$$
On the other hand, it follows from \eqref{edual} that for $f\in S(V;\A)\subset V\ltimes A$ the twisted dual action is defined by
$$
\hat\alpha^{\Omega_J}_y(f)(x)=e(-x\cdot y)\alpha_{-Jy}(f(x)),
$$
since $r_{\Omega_J}=J$. It is routine to check that $\Theta_J$ intertwines these two actions.
\ep

\begin{example}
Assume $A=C_0(V)$ and $\alpha$ is the action by left translations, so $\alpha_x(f)(y)=f(y-x)$. In this case the previous theorem gives the well-known isomorphism $C_0(V)_J\cong C^*_r(V;\Omega_J)$. Indeed, we have
$$
C_0(V)_{\Omega_J}=E(C^*_r(V;\Omega_J)\otimes1)E^*\cong C^*_r(V;\Omega_J),
$$
where we consider $C^*_r(V;\Omega_J)\otimes C_0(V)$ as an algebra of operators on $L^2(V)\otimes L^2(V)$ and define a unitary $E$ on this Hilbert space by $(E\xi)(x,y)=e(x\cdot y)\xi(x,y)$. A simple computation shows that if $f\in\A$ lies in the Fourier algebra of $V$, then
$$
\Ad(E^*\F\otimes1)\pi_J(f)=\int_V\hat f(x)(\lambda^\Omega_x\otimes1)dx,
$$
so the isomorphism $C_0(V)_J\cong C^*_r(V;\Omega_J)$ we thus obtain is defined by the familiar formula $f\mapsto \int_V\hat f(x)\lambda^{\Omega_J}_xdx$.
\ee
\end{example}

Theorem~\ref{tdual} provides a short proof of equivalence of the approaches of Rieffel and Kasprzak to deformation. There are, however, a couple of loose ends left to tie. In \cite[Theorem~A.3]{BNS} we already constructed an isomorphism between $A_J$ and $A_{\Omega_J}$. A natural question is whether this is the same isomorphism. A related, and more interesting, question is how the quantization maps $(T_\nu\otimes\iota)\alpha \colon A\to A_{\Omega_J}$ look like in Rieffel's picture.

Consider a normal linear functional $\nu$ on $W^*(V;\bar\Omega_J)$. Assume that the function $x\mapsto\nu(\lambda^{\bar\Omega_J}_x)$ lies in the Fourier algebra of $V$, so it is the Fourier transform of a function $g_\nu\in L^1(V)$. We then define a linear map
$$
\Phi_\nu\colon \A\to\A_J\ \ \text{by}\ \ \Phi_\nu(a)=\int_V\alpha_x(a)g_\nu(x)dx.
$$
Note that if $\nu=(\cdot\,\xi,\zeta)$ for some $\xi,\zeta\in L^2(V)$, then it is easy to check that the assumption on $\nu$ is satisfied if e.g.~$\hat\xi,\zeta\in L^1(V)$, and we have
$$
g_\nu(x)=\int_V\hat\xi(x+Jy)\overline{\zeta(y)}e(x\cdot y)dy,\ \ \nu(\lambda^{\bar\Omega_J}_x)=\int_V\xi(y-x)\overline{\zeta(y)}e(x\cdot Jy)dy.
$$

\begin{proposition} \label{pintert}
Assume $\nu\in W^*(V;\bar\Omega_J)_*$ is such that the function $x\mapsto\nu(\lambda^{\bar\Omega_J}_x)$ lies in the Fourier algebra of $V$. Then for any $a\in\A$ we have
\begin{equation} \label{eintert}
\pi_J(\Phi_\nu(a))=\Ad(\F^*\otimes1)(T_\nu\otimes\iota)\alpha(a).
\end{equation}
\end{proposition}

\bp We may assume that $A\subset B(H)$ for some Hilbert space $H$ and that the action $\alpha$ of $V$ on $A$ is implemented by a strongly continuous unitary representation of $V$ on $H$. Then $\pi_J$ can be considered as a representation of $A_J$ on $L^2(V)\otimes H=L^2(V;H)$, so that we have
\begin{equation} \label{erep}
(\pi_J(a)\xi)(x)=\int_V\alpha_{Jy-x}(a)\hat\xi(y)e(x\cdot y)dy\ \ \text{for}\ \ a\in\A_J\ \ \text{and}\ \ \xi\in S(V;H).
\end{equation}
The right hand side of \eqref{eintert} is a normal map in $a\in\A\subset B(H)$. On the other hand, if $\{a_i\}_i$ is a bounded net in $\A$ converging strongly to $a\in\A$, then it follows from \eqref{erep} that $(\pi_J(\Phi_\nu(a_i))\xi,\zeta)\to (\pi_J(\Phi_\nu(a))\xi,\zeta)$ for all $\xi,\zeta\in S(V;H)$. Therefore it suffices to check \eqref{eintert} for elements of a strongly dense subset of the unit ball of $\A$.

Using the above observation and the same trick based on the Takesaki duality as in the proof of \cite[Lemma~A.1]{BNS}, it is now easy to prove the proposition. Namely, by embedding $A$ into a much larger algebra we may assume that the spectral subspaces of $A$ span a strongly dense subalgebra of~$A$. Hence it suffices to check~\eqref{eintert} on homogeneous elements. If $a\in A_z$, so that $\alpha_x(a)=e(-x\cdot z)a$, then by the computation in Example~\ref{exap} we have
$$
\Ad(\F^*\otimes1)(T_\nu\otimes\iota)\alpha(a)=\nu(\lambda^{\bar\Omega}_z)\alpha(a)(\lambda_{-Jz}\otimes1)
=\hat g_\nu(z)\alpha(a)(\lambda_{-Jz}\otimes1).
$$
On the other hand, $\Phi_\nu(a)=\hat g_\nu(z)a$ and
$$
(\pi_J(a)\xi)(x)=\int_Va\hat\xi(y)e(x\cdot y)e((x-Jy)\cdot z)dy=\alpha_{-x}(a)\xi(x+Jz)\ \ \text{for}\ \ \xi\in S(V;H).
$$
But this is exactly how the operator $\alpha(a)(\lambda_{-Jz}\otimes1)$ acts.
\ep

In our current notation Theorem~A.3 in~\cite{BNS} can be formulated by saying that if $J^2=-\pi^2h^2$ for some number $h>0$, then there exists a unique isomorphism $A_J\cong A_{\Omega_J}$ such that
\begin{equation} \label{eiso}
\Phi(a)=\frac{1}{(\pi h)^{d/2}}\int_Ve^{-\frac{1}{h}\|x\|^2}\alpha_x(a)dx\mapsto (T_{\nu_0}\otimes\iota)\alpha(a)\ \ \text{for all}\ \ a\in A,
\end{equation}
where $\nu_0$ is the normal state on $W^*(V;\bar\Omega_J)$ defined by $\nu_0(\lambda^{\bar\Omega_J}_x)=e^{-\pi^2h\|x\|^2}$. Since
$$
g_{\nu_0}(x)=\frac{1}{(\pi h)^{d/2}}e^{-\frac{1}{h}\|x\|^2},
$$
we have $\Phi=\Phi_{\nu_0}$, so the map \eqref{eiso} coincides with the isomorphism from Theorem~\ref{tdual}. As was already remarked in \cite{BNS}, the state $\nu_0$ is the vacuum state on the algebra of canonical commutation relations. Therefore the map $(T_{\nu_0}\otimes\iota)\alpha\colon A\to A_{\Omega_J}$ is the most natural among the quantization maps $(T_{\nu}\otimes\iota)\alpha$. At the same time we see now that the vacuum state, as well as the map $\Phi$ introduced in~\cite{KNW}, does not play any special role in constructing the isomorphism $A_J\cong A_{\Omega_J}$. In particular, we have the following immediate corollary to Proposition~\ref{pintert}, which extends a result in \cite{KNW} for the map~$\Phi$ and is valid for any skew-symmetric $J$.

\begin{corollary}
For any $\nu\in W^*(V;\bar\Omega_J)_*$ such that the function $x\mapsto\nu(\lambda^{\bar\Omega_J}_x)$ lies in the Fourier algebra of $V$, the map $\Phi_\nu\colon \A\to\A_J$ extends to a completely bounded map $A\to A_J$, and $\|\Phi_\nu\|_{\rm cb}\le\|\nu\|$. If in addition $\nu$ is positive, then $\Phi_\nu$ is completely positive.
\end{corollary}

\bigskip

\end{document}